\documentclass[11pt]{amsart}
\usepackage{amsmath}

\numberwithin{equation}{section}

\theoremstyle{plain}

\newtheorem{Main}{Main Theorem}

\newtheorem{MainCor}{Corollary}

\newtheorem{Prop}{Proposition}[section]

\newtheorem{Lem}[Prop]{Lemma}

\theoremstyle{remark}

\newcommand{\CBbb}{\Bbb C}
\newcommand{\RBbb}{\Bbb R}

\newcommand{\End}{{\rm End}\, }
\newcommand{\ad}{{\rm ad}\, }

\newcommand{\lra}{\longrightarrow}

\newcommand{\Hol}{\mathcal{H}}
\newcommand{\A}{\mathcal{A}}
\newcommand{\Oka}{\mathcal{O}}
\newcommand{\G}{\mathfrak{G}}
\newcommand{\Gc}{{\G}^{\CBbb}}
\newcommand{\B}{\mathfrak{B}_\tau}
\newcommand{\M}{\mathfrak{M}(r,d)}
\newcommand{\bn}{\rho^{k-1}_{r,d}}
\newcommand{\bnlocus}{W^{k-1}_{r,d}}

\begin{document}


\title[Brill-Noether Problem for Vector Bundles]
	{On the Brill-Noether Problem \\ for Vector Bundles}

\author[G. D. Daskalopoulos]{Georgios D. Daskalopoulos}

\address{Department of Mathematics \\
		Brown University \\
		Providence,  RI  02912}

\email{daskal@gauss.math.brown.edu}

\thanks{Supported in part by NSF grant DMS-9504297}

\author[R. A. Wentworth]{Richard A. Wentworth}

\address{Department of Mathematics \\
   University of California \\
   Irvine,  CA  92717}

\email{raw@math.uci.edu}

\thanks{Supported in part by NSF grant DMS-9503635 and a Sloan Fellowship}

\date{October 4, 1995}


\begin{abstract}
On an arbitrary compact Riemann surface,
necessary and sufficient conditions are found for the existence of semistable
vector
bundles  with slope between zero and one and a prescribed number of
linearly independent holomorphic sections.  Existence is achieved by minimizing
the
Yang-Mills-Higgs functional.
\end{abstract}

\maketitle

\section{Introduction}

In this note we exhibit the existence of semistable vector bundles on compact
Riemann
surfaces with a prescribed number of linearly independent holomorphic sections.
 This may be
regarded as a higher rank version of the classical Brill-Noether problem for
line bundles.

Fix a compact Riemann surface $\Sigma$ of genus $g\geq 2$ and integers $r$ and
$d$ satisfying
\begin{equation}
0\leq d \leq r\ , \quad r\geq 2\ .  \label{d-range}
\end{equation}
Then the main result may be stated as follows:
\begin{Main}
Let $k$ be a positive integer and suppose that $r$ and $d$ satisfy
\eqref{d-range}.  Then
the necessary and sufficient conditions for the existence of a semistable
bundle on
$\Sigma$ with at least $k$ linearly independent holomorphic sections are $k\leq
r$ and if
$d\neq 0$,
$r\leq d+(r-k)g$.
\end{Main}

By analogy with the classical situation of special divisors (cf.\
\cite{ACGH,N}) one can
define the higher rank version of the Brill-Noether number:
\begin{equation}
\bn= r^2(g-1)+1-k(k-d+r(g-1))\ .
\end{equation}
Then $\bn$ is the formal dimension of the locus $\bnlocus$ in the moduli space
of
semistable bundles of rank $r$ and degree $d$.  $\bnlocus$ is defined as the
closure of the
set of stable bundles with at least $k$ linearly independent sections.  Note
that the
condition in the Main Theorem implies that $\bn\geq 1$, except in the trivial
case $d=0$
where $\bnlocus$ is necessarily empty.  The converse, in general, is not true.
Thus, unlike
the case of divisors, there are situations where
$\bn\geq 0$ and
$\bnlocus=\emptyset$.

It would be interesting to improve the Main Theorem to a statement concerning
stable
bundles; however, our method does not immediately imply such a result except in
 special
cases. We do have the following:
\begin{MainCor} {\rm (i) (see \cite[Thm III.2.4]{S})}  For  $d>0$ and any rank
$r$, there
exists a stable bundle of rank $r$ and degree $d$ with a nontrivial holomorphic
section.
{\rm (ii)}   If $0<d<r$  and $r\leq d+g$, then there exists a stable bundle of
rank
$r$ and degree $d$ with precisely $r-1$ linearly independent sections.
\end{MainCor}

Instead of the constructive approach to theorems of this type taken in
references
\cite{S,T}, we use a variational method.  More precisely, we study the Morse
theory of the
Yang-Mills-Higgs functional (cf.\ \cite{B}).  The idea is simply the following:
 Let
$(A^i,\vec\varphi^i)$ be a minimizing sequence with respect to the
Yang-Mill-Higgs
functional~\eqref{YMH}.  Here,
$\vec\varphi_i=(\varphi^i_1,\ldots,\varphi^i_k)$ is a $k$-tuple of linearly
independent
holomorphic sections with respect to $A^i$.  The minimal critical values
correspond to
solutions to the $k$-$\tau$-vortex equations, which for an appropriate choice
of $\tau$ imply that the limiting holomorphic structure is semistable (cf.\
\cite{BDW}).
If the sequence is assumed to converge to a nonminimal critical value, then we
show that
under the assumptions  of the Main Theorem there exist ``negative directions"
which contradict the fact that the sequence is minimizing.

The energy estimates used closely follow  \cite{D}.  However, an extra
combinatorial
argument is needed to ensure that the bundles constructed have the correct
number of
holomorphic sections, and this is where the assumption $r\leq d+(r-k)g$ is
needed.

We have been informed that the Main Theorem stated above has been proven using
somewhat different methods in \cite{BGN}.

\medskip
\noindent {\it Acknowledgements.}  The authors would like to thank L. Brambila
Paz for
introducing them to this problem and for several useful discussions during the
preparation of
this manuscript.  They are also grateful for the warm hospitality of UAM,
Mexico and the
Max-Planck Institute in Bonn, where a portion of this work was completed.

\section{The Yang-Mills-Higgs Functional}       \label{S:functional}

Let $\Sigma$, $d$, and $r$ be as in the Introduction, and let $k$ be a positive
integer.
Let $E$ be a fixed hermitian vector bundle on $\Sigma$ of rank $r$ and degree
$d$.  Let $\A$
denote the space of hermitian connections on $E$, $\Omega^0(E)$ the space of
smooth sections
of $E$, and $\Hol\subset \A\times \Omega^0(E)^{\oplus k}$ the subspace
consisting of
holomorphic $k$-pairs.  Thus,
$$
\Hol=\left\{ \left(A, \vec\varphi=(\varphi_1,\ldots,\varphi_k)\right) :
 D''_A\varphi_i=0 , \
i=1,\ldots,k\right\}\ .
$$
Given a real parameter $\tau$, we define the Yang-Mills-Higgs functional:
\begin{align}
f_\tau &: \A\times\Omega^0(E)^{\oplus k}\longrightarrow \RBbb \notag \\
f_\tau(A,\vec\varphi)&=
\Vert F_A\Vert^2 + \sum_{i=1}^k\Vert D_A\varphi_i\Vert^2 +\frac{1}{4}\left\Vert
\sum_{i=1}^k \varphi_i\otimes\varphi_i^\ast-\tau {\bf I}\, \right\Vert^2
-2\pi\tau d
\label{YMH}
\end{align}
In the above, the $\Vert\cdot\Vert$ denotes $L^2$ norms.  Using a Weitzenb\"ock
formula
we obtain (cf.\ \cite[Theorem 4.2]{B})
$$
f_\tau(A,\vec\varphi)=2\sum_{i=1}^k\Vert D''_A
\varphi_i\Vert^2+\left\Vert\sqrt{-1}\Lambda
F_A+\frac{1}{2}\sum_{i=1}^k\varphi_i\otimes\varphi_i^\ast-\frac{\tau}{2}{\bf
I}\,
\right\Vert^2\ ,
$$
and therefore the absolute minimum of $f_\tau$ consists of holomorphic
$k$-pairs satisfying
the $k$-$\tau$-vortex equations discussed in \cite{BDW}.
\begin{Prop} \label{P:gradient}
{\rm (i)}  The $L^2$-gradient of $f_\tau$ is given by
\begin{align*}
\left(\nabla_{(A,\vec\varphi)} f_\tau\right)_1 &=
D_A^\ast F_A+\frac{1}{2}\sum_{j=1}^k\left(
D_A\varphi_j\otimes\varphi_j^\ast-\varphi_j\otimes D_A\varphi_j^\ast\right) \\
\left(\nabla_{(A,\vec\varphi)} f_\tau\right)_{2,i} &=
\Delta_A\varphi_i-\frac{\tau}{2}\varphi_i+\frac{1}{2}\sum_{j=1}^k\langle
\varphi_i,\varphi_j\rangle\varphi_j
\end{align*}
{\rm (ii)}  If $(A,\vec\varphi)\in\Hol$, then under the usual identification
$\Omega^0(\Sigma,\ad E)\simeq$ \break $\Omega^{0,1}(\Sigma,\End E)$,
we have
\begin{align*}
\left(\nabla_{(A,\vec\varphi)} f_\tau\right)_1 &=
 -D''_A\left(\sqrt{-1}\Lambda
F_A+\frac{1}{2}\sum_{j=1}^k\varphi_j\otimes\varphi_j^\ast\right) \\
\left(\nabla_{(A,\vec\varphi)} f_\tau\right)_{2,i} &=
\sqrt{-1}\Lambda
F_A(\varphi_i)-\frac{\tau}{2}\varphi_i+\frac{1}{2}\sum_{j=1}^k\langle
\varphi_i,\varphi_j\rangle\varphi_j
\end{align*}
{\rm (iii)}  If $(A,\vec\varphi)\in\Hol$ is a critical point of $f_\tau$, then
either {\rm (I)}
$\vec\varphi\equiv 0$ and $A$ is a direct sum of Hermitian-Yang-Mills
connections (not
necessarily of the same slope), or {\rm (II)} $A$ splits as $A=A'\oplus A_Q$ on
$E=E'\oplus E_Q$, where $(A',\vec\varphi)$ solves the $k$-$\tau$-vortex
equations and $A_Q$
is a direct sum of Hermitian-Yang-Mills connections (not necessarily of the
same slope).
\end{Prop}

\begin{proof}
(i) is a standard calculation, and (ii) follows from (i) via the K\"ahler
identities.  We
are going to prove (iii).  If $(A,\vec\varphi)$ is critical, then since
$
\sqrt{-1}\Lambda F_A+\frac{1}{2}\sum_{j=1}^k\varphi_j\otimes\varphi_j^\ast
$
is a self-adjoint holomorphic endomorphism, it must give a splitting
$A=A_0\oplus\cdots\oplus A_\ell$ according to its distinct (constant)
eigenvalues
$\sigma_0,\ldots,\sigma_\ell$.  Write
$$
\sqrt{-1}\Lambda F_A=
\begin{pmatrix}
-\frac{1}{2}\sum_{j=1}^k\varphi_j\otimes\varphi_j^\ast +\sigma_0\ {\bf I}
& 0 & \cdots & 0 \\
0& \sigma_1\ {\bf I} & & \vdots \\
\vdots && \ddots & \\
0 & \cdots && \sigma_\ell\ {\bf I}
\end{pmatrix}
\quad .
$$
Thus,
\begin{align*}
0 &=
\sqrt{-1}\Lambda
F_A(\varphi_i)-\frac{\tau}{2}\varphi_i+\frac{1}{2}\sum_{j=1}^k\langle
\varphi_i,\varphi_j\rangle\varphi_j\\
&= -\frac{1}{2}\sum_{j=1}^k\langle\varphi_i,
\varphi_j\rangle\varphi_j+\sigma_0\varphi_i
-\frac{\tau}{2}\varphi_i+\frac{1}{2}\sum_{j=1}^k\langle
\varphi_i,\varphi_j\rangle\varphi_j\\
&= \left(\sigma_0-\frac{\tau}{2}\right)\varphi_i\ ,
\end{align*}
from which we obtain either Case I or Case II, depending upon whether
$\vec\varphi\equiv 0$.
\end{proof}

Next, recall that $\Hol$ is an infinite dimensional complex analytic
variety whose tangent space is given by the kernel of a certain differential
defined in
\cite[3.15]{BDW}.  Moreover, $\Hol$  is preserved by the action of the complex
gauge
group $\Gc$.  We have the following:
\begin{Prop}
\label{P:tangent}
If $(A,\vec\varphi)\in \Hol$, then $\nabla_{(A,\vec\varphi)}f_\tau$ is tangent
to the
orbits of $\Gc$.  In particular, $\nabla_{(A,\vec\varphi)}f_\tau$ is tangent to
$\Hol$
itself.
\end{Prop}

\begin{proof}
Set
$
u=\sqrt{-1}\Lambda
F_A+\frac{1}{2}\sum_{j=1}^k\varphi_j\otimes\varphi_j^\ast-\frac{\tau}{2}{\bf I}
$.
By Proposition \ref{P:gradient} (ii) we have that
$\nabla_{(A,\vec\varphi)}f_\tau=d_1(u)$,
where
$d_1$ is the differential defined in \cite[3.15]{BDW}.  The Proposition
follows.
\end{proof}

Because of Proposition \ref{P:tangent}, the critical points of the functional
$f_\tau$
restricted to $\Hol$ are characterized by Proposition \ref{P:gradient} (iii).

A solution $(A(t),\vec\varphi(t))$, $t\in [0,t_0)$ to the initial value problem
\begin{equation}
\label{flow}
\left(\frac{\partial A}{\partial t}, \frac{\partial\vec\varphi}{\partial
t}\right)
=
-\nabla_{(A,\vec\varphi)}f_\tau\ , \quad
\left( A(0),\vec\varphi(0)\right)
=
(A_0,\vec\varphi_0)\ ,
\end{equation}
is called the $L^2$-gradient flow of $f_\tau$ starting at
$(A_0,\vec\varphi_0)$.  Notice
that
\begin{equation}
\label{E:decay}
\frac{d}{dt}f_\tau(A(t),\vec\varphi(t))
=-\left\Vert\nabla_{(A(t),\vec\varphi(t))}
f_\tau\right\Vert^2\ ,
\end{equation}
and so the energy decreases along the $L^2$-gradient flow.
\begin{Prop}
\label{P:shorttime}
Given $(A_0,\vec\varphi_0)\in\Hol$, there is a $t_0>0$ such that the
$L^2$-gradient flow
exists for $0\leq t< t_0$.
\end{Prop}

\begin{proof}
The proof is an application of the implicit function theorem as in \cite{R}.
\end{proof}

\section{Technical Lemmas}
\label{S:technical}

In this section we collect several results needed for the proof of the Main
Theorem.
Throughout,
$E$ will denote a holomorphic bundle of rank $r$ and degree $d$ on the compact
Riemann
surface $\Sigma$.

\begin{Lem}
\label{L:bounds}
Let $E$ be as above with $0\leq d\leq r$  and $h^0(E)=k$.  If either {\rm (i)}
$E$ is semistable, or {\rm (ii)} $E$ satisfies the k-$\tau$-vortex equation
for some $0<\tau<1$ and
$E$ does not contain the trivial bundle as a split factor; then $k\leq  r$ and
if $d\neq 0$, $r\leq d+(r-k)g$.
\end{Lem}

\begin{proof}  We first show that $k\leq r$.  Suppose $k\geq r$. Thus, $E$ has
at least $r$
linearly independent holomorphic sections. If the sections generate
$E$ at every point, then $E\simeq\Oka^{\oplus r}$; in which case $d=0$ and
$k=r$.
Suppose the sections fail to generate at every point.  Then we
can find a point
$p\in\Sigma$ and a section of $E$ vanishing at $p$.  Thus $E$ contains
$\Oka(p)$ as a
subsheaf, which is a contradiction to (ii).  If (i) is assumed, then $E$ is
strictly
semistable with $d=r$, and the bound $k\leq r$ follows from induction on the
rank. Note
that the second inequality is also satisfied in this case.

Assume $0< d< r$.  In both cases (i) and (ii) we obtain
$0\to \Oka^{\oplus k} \xrightarrow{\pi} E\to F\to 0$,
 where $F$ is locally
free.  By dualizing and taking the resulting long exact sequence in cohomology,
we find
$$
0\lra H^0(F^\ast)\lra H^0(E^\ast)
\xrightarrow{\delta} H^0(\Oka^{\oplus k}) \lra H^1(F^\ast)
\ .
$$
We are going to show that $H^0(E^\ast)=0$.  The result then follows by the
Riemann-Roch
formula.  For
(i), $H^0(E^\ast)=0$ by semistability.  For (ii), note first that $\delta=0$.
For if not,
there would be a section $s:\Oka\to E^\ast$ with $\pi^\ast\circ s=\sigma\neq
0$.  But
$\sigma$ could not have any zeros, and so $\Oka$ would  be a split factor in
$E^\ast$; hence,
also in $E$. Secondly, we show that $H^0(F^\ast)=0$. Let
$L\subset F^\ast$ be a subbundle.  Then  $\tau$-stability immediately implies
$c_1(L^\ast) >
\tau >0$.  Thus, in particular, $F^\ast$ cannot contain $\Oka$ as  a subsheaf.
This
completes the  proof.
\end{proof}

\begin{Lem}
\label{L:extension}
Let $E_1$, $E_2$ be holomorphic bundles of rank $r_1, r_2$ and degree $d_1,
d_2$,
satisfying $0\leq \mu_1=d_1/r_1 < d_2/r_2=\mu_2\leq 1$.  Suppose
$h^0(E_1)=k_1 \leq r_1$, $h^0(E_2)=k_2\leq r_2$,  and
$$
d_2+(r_2-k_2-1)g < r_2\leq d_2+(r_2-k_2)g\ .
$$
Furthermore,
\begin{itemize}
\item  If $d_1\neq 0$  assume $r_1\leq d_1+(r_1-k_1)g$.
\item  If $d_1=0$ and $k_1=r_1$, assume $r_2 < d_2+(r_2-k_2)g$.
\end{itemize}
Then there exists a nontrivial extension $0\to E_1\to E\to E_2\to 0$ such that
$h^0(E)=
k_1+k_2$.
\end{Lem}

\begin{proof} If $k_2=0$, the result follows from Riemann-Roch.  Suppose
$k_2\geq 1$. The
condition that the
$k_2$ sections of
$E_2$ lift for some nontrivial extension is
$k_2 h^1(E_1) < h^1(E_1\otimes E_2^\ast)$.  Notice that
\begin{align*}
h^1(E_1) &= h^0(E_1)-d_1+r_1(g-1)=k_1-d_1+r_1(g-1) \\
h^1(E_1\otimes E_2^\ast) &= h^0(E_1\otimes E_2^\ast)+r_1 r_2(\mu_2-\mu_1+g-1)
\\
&\geq  r_1 r_2(\mu_2-\mu_1+g-1)\ ,
\end{align*}
hence, it suffices to show that
$$
k_2(k_1-d_1+r_1(g-1)) < r_1 r_2(\mu_2-\mu_1+g-1)\ ,
$$
or equivalently, that
\begin{equation}
\label{E:one}
r_1(d_2-r_2+(r_2-k_2)g)-r_2 d_1 + k_2 d_1 -k_1 k_2 + k_2 r_1 > 0\ .
\end{equation}
Now if $k_2=r_2=d_2$, then \eqref{E:one} is trivially satisfied by the
hypotheses.
Similarly for $d_1=0$. So assume
$k_2\leq r_2-1$, $d_1\neq 0$. Write $d_2=r_2-(r_2-k_2)g+p$, where $0\leq p < g$
by
assumption.  On the other hand,
$$
d_1 < r_1\frac{d_2}{r_2}\leq
\left(d_1+(r_1-k_1)g\right)\frac{r_2-(r_2-k_2)g+p}{r_2}\ ,
$$
which is equivalent to
$$
-\frac{d_1 p}{g}+k_1 p + (r_1-k_1)(r_2-k_2)(g-1) < r_1 p-r_2d_2+k_2 d_1-k_1
k_2+k_2 r_1
\ .
$$
Therefore, \eqref{E:one} will follow from
\begin{equation}
\label{E:two}
-\frac{d_1 p}{g}+k_1 p + (r_1-k_1)(r_2-k_2)(g-1)\geq 0\ .
\end{equation}
Now if $p=0$ then \eqref{E:two} is trivially satisfied.  Assume that $1\leq
p\leq g-1$.
Then
\begin{align*}
-\frac{d_1 p}{g} &+k_1 p + (r_1-k_1)(r_2-k_2)(g-1) \\
&\geq -d_1+r_1
p-(r_1-k_1)p+(r_1-k_1)(r_2-k_2)(g-1) \\
&\geq (r_1-d_1)+(r_1-k_1)(r_2-k_2-1)(g-1)  \\
&\geq 0 \ ,
\end{align*}
which proves \eqref{E:two}, \eqref{E:one}, and hence the Lemma.
\end{proof}

In order to get an upper bound on the infimum of the Yang-Mills-Higgs
functional in the next
section, we shall need the following construction and energy estimate:
\begin{Lem}
\label{L:special}
Assume $0< d < r$, $k\geq 1$, and $r\leq d+(r-k)g$.  Let $F$ be a holomorphic
bundle of
degree
$d$ and rank
$r-1$ with $h^0(F)=k-1$.  Then there exists a non-split extension $0\to\Oka\to
E\to F\to 0$
with $h^0(E)=k$.
\end{Lem}

\begin{proof}
The condition for all of the sections of $F$ to lift is
\begin{align*}
(k-1)h^1(\Oka) <  h^1(F^\ast)\
&\iff\quad g(k-1) < d+(r-1)(g-1) \\
&\iff\quad r< d+(r-k)g+1\ ,
\end{align*}
and hence the result.
\end{proof}

\begin{Prop}[{cf.\ \cite[Prop. 3.5]{D}}]
\label{P:energyestimate}
Let $E_1, E_2$ be hermitian bundles with slope $\mu_1, \mu_2$.  Let $A_1, A_2$
be hermitian
connections on $E_1, E_2$, and $\vec\varphi^1, \vec\varphi^2$ be $k_1$ and
$k_2$ tuples of
holomorphic sections.  Set $k=k_1+k_2$. Let  $\tau_1, \tau_2$ and $\tau$ be
real numbers
satisfying $\mu_1\leq\tau_1\leq \tau < \mu_2 \leq\tau_2$, and assume that
$(A_1,\vec\varphi^1)$ and
$(A_2,\vec\varphi^2)$ satisfy the $\tau_1$ and $\tau_2$ vortex equations,
respectively.
Set $E=E_1\oplus E_2$, $\varphi_i=(\varphi^1_i,0)$ for $i=1,\ldots, k_1$, and
$\varphi_{k_1+i}=(0,\varphi^2_i)$ for $i=1,\ldots, k_2$.  Then there exist
constants
$\varepsilon_1, \varepsilon_2 , \eta >0$ such that for all
$$
\beta\in  H^{0,1}\left(\Sigma, {\rm Hom}(E_2,E_1)\right)\ ,\quad
\vec\psi\in\Omega^0(E)^{\oplus k}\ ,
$$
with $\Vert\beta\Vert =\varepsilon_1$, $\Vert\vec\psi\Vert\leq \varepsilon_2$,
and
$$
\left( A_\beta=
\begin{pmatrix}
A_1 & \beta \\ 0 & A_2
\end{pmatrix},
\vec\varphi+\vec\psi\right)\in\Hol\ ,
$$
it follows that
$
f_\tau(A_\beta,\vec\varphi+\vec\psi) <  f_\tau(A_1\oplus A_2,\vec\varphi)-\eta
$.
\end{Prop}

\begin{proof}
By assumption,
$$
\sqrt{-1}\Lambda
F_{A_\ell}+\frac{1}{2}\sum_{j=1}^{k_\ell}
\varphi^{\ell}_j\otimes(\varphi^{\ell}_j)^\ast
=
\frac{\tau}{2}{\bf I}_{\ell}\ ,\quad \ell=1,2\ .
$$
It follows that
$$
\sqrt{-1}\Lambda
F_{A_1\oplus
A_2}+\frac{1}{2}\sum_{j=1}^{k_1}\varphi^1_j\otimes(\varphi^1_j)^\ast
+\frac{1}{2}\sum_{j=1}^{k_2}\varphi^2_j
\otimes(\varphi^2_j)^\ast-\frac{\tau}{2}{\bf I}=
\begin{pmatrix}
\frac{\tau_1-\tau}{2}{\bf I}_1 & 0 \\
0 & \frac{\tau_2-\tau}{2}{\bf I}_2
\end{pmatrix}
$$
The argument of \cite[pp.\ 715-716]{D} shows that there is a constant
$\varepsilon_1$  such
that for $\beta$ and $A_\beta$ as in the statement,
$$
f_\tau\left(
A_\beta,\varphi^1_1,\ldots,\varphi^1_{k_1},
\varphi^2_1,\ldots,\varphi^2_{k_2}\right)
<
f_\tau\left(
A_1\oplus
A_2,\varphi^1_1,\ldots,\varphi^1_{k_1},
\varphi^2_1,\ldots,\varphi^2_{k_2}\right)
\ .
$$
Now if we take $\varepsilon_2$ sufficiently small the Proposition follows (note
that which
norms we use is irrelevent, since $\beta$ and $\vec\varphi+\vec\psi$ satisfy
elliptic
equations, and hence the $L^2$ norm is equivalent to any other).
\end{proof}

\section{Proof of the Main Theorem}
\label{S:proof}

Necessity of the conditions follows from Lemma \ref{L:bounds}, and sufficiency
for $d=0$ or
$d=r$ is clear as well. To prove sufficiency in general, we shall proceed by
induction on the
rank.  The case
$r=2$,
$d=1$ is clear from a direct construction.  Assume the Main Theorem holds for
bundles of
rank $< r$.  We show that it holds for $r$ as well.
Let $\Hol^\ast\subset\Hol$ denote the subspace of $k$-pairs
$\left(A,\vec\varphi=(\varphi_1,\ldots,\varphi_k)\right)$
such that the sections $\varphi_1,\ldots,\varphi_k$ are linearly independent.
Fix $\tau$ as
in Assumption 1 of \cite{BD}, i.e. $\mu(E)<\tau=\mu(E)+\gamma < \mu_+$, where
$\mu_+$ is
the smallest possible slope greater that $\mu=\mu(E)$  of a subbundle of $E$
(note that
$0<\tau<1$ and that we also normalize the volume of $\Sigma$ to be $4\pi$).

\begin{Lem}
\label{L:inf}
Let $m=\inf_{\Hol^\ast} f_\tau$.  Then $0\leq m<\pi/(r-1)$.
\end{Lem}

\begin{proof}
Let $F$ be a vector bundle of degree $d$ and rank $r-1$. Then by the
inductive hypothesis, we may assume there exist hermitian connections $A_1$ and
$A_2$ on
$\Oka$ and
$F$, respectively, and holomorphic sections $\varphi_1\neq 0$ in
$H^0(\Sigma,\Oka)$, and
$\varphi_2,\ldots,\varphi_k$ linearly independent sections in $H^0(\Sigma,F)$,
such that
$(A_1,\varphi_1)$ and $(A_2,\varphi_2,\ldots,\varphi_k)$ satisfy the $\tau_1$
and
$\tau_2$ vortex equations, respectively, for $\tau_1=\tau$,
$\tau_2=d/(r-1)+\gamma$.
It follows from Lemma \ref{L:special} and Proposition \ref{P:energyestimate}
that there is
a nontrivial extension $\beta: 0\to\Oka\to E\to F\to 0$, and $\vec\psi$ such
that
$(A_\beta,\vec\varphi+\vec\psi)\in\Hol^\ast$ and
\begin{align*}
f_\tau(A_\beta,\vec\varphi+\vec\psi)&<f_\tau(A_1\oplus
A_2,\varphi_1,\ldots,\varphi_k)-\eta \\
&=
\left\Vert\frac{1}{2}\left(\frac{d}{r-1}-\frac{d}{r}\right){\bf I_F}\,
\right\Vert^2-\eta <
\frac{\pi}{r-1}\ .
\end{align*}
\end{proof}

Let $(A^i,\vec\varphi^i)$ be a minimizing sequence in $\Hol^\ast$.  Thus,
$f_\tau(A^i,\vec\varphi^i)\to m$.  By weak compactness (more precisely, see the
argument
in \cite[Lemma 5]{BD}) there is a subsequence converging to
$(A^\infty,\vec\varphi^\infty)$ in the $C^\infty$ topology.  By the continuity
of
$f_\tau$ with respect to the $C^\infty$ topology, Propositions
\ref{P:shorttime} and
\ref{P:tangent}, and equation \eqref{E:decay}, it follows that
$(A^\infty,\vec\varphi^\infty)$ is a critical point of
$f_\tau$.  If  the holomorphic structure $E^\infty$ defined by $A^\infty$ is
semistable, then
by semicontinuity of cohomology  we are finished.  We therefore
assume
$E^\infty$ is unstable and derive a contradiction.  According to Proposition
\ref{P:gradient} (iii) we must consider the following cases:
\begin{align*}
\vec\varphi^\infty &= 0\ , \quad E^\infty=E_1\oplus\cdots\oplus E_\ell
\tag{I}   \\
\vec\varphi^\infty &\neq 0\ , \quad E^\infty=E_{\varphi}\oplus
E_1\oplus\cdots\oplus E_\ell
\tag{II}
\end{align*}
Set $\mu_j=\mu(E_j)$, and assume $\mu_1 < \cdots < \mu_\ell$.  If $\mu_\ell >
1$ (or
similarly,
$\mu_1 < 0$), then
$$
f_\tau\left(A^\infty,\vec\varphi^\infty\right)\geq \pi(\mu_\ell-\tau)^2 r_\ell
\geq \pi(\mu_\ell-1)^2 r_\ell
\geq\frac{\pi}{r_\ell}\geq\frac{\pi}{r-1} > m\ ,
$$
contradicting Lemma \ref{L:inf}.
We therefore rule out this possibility.
We will consider Cases I and II separately.

\medskip
\noindent \emph{Case I}\@.  Let $k_i=h^0(E_i)$.  By semicontinuity of
cohomology,
$\sum_{i=1}^\ell k_i\geq k$.  If $\mu_\ell =1$, then we may replace $E_\ell$ by
a
Hermitian-Yang-Mills bundle $\widehat E_\ell$ with exactly $\hat k_\ell=r_\ell$
sections.
Hence, we may assume that
$$
d_\ell + (r_\ell -\hat k_\ell-1)g < r_\ell \leq d_\ell +(r_\ell-\hat k_\ell)g\
{}.
$$
For $1<i<\ell$, the inductive hypothesis implies that we may replace $E_i$ by a
Hermitian-Yang-Mills bundle $\widehat E_i$ with
$$
h^0(\widehat E_i)=\hat k_i=\left[ \frac{d_i+r_i(g-1)}{g}\right]\ ,
$$
the maximal number of sections allowed for $d_i, r_i$, and $g$.
Note that
\begin{equation}
\label{E:max}
d_i + (r_i -\hat k_i-1)g < r_i \leq d_i +(r_i-\hat k_i)g\ .
\end{equation}
 If $\mu_1\neq 0$, then we
can replace $E_1$ by $\widehat E_1$ as above.  If $\mu_1=0$, we may replace
$E_1$ with
$\Oka^{\oplus r_1}$, with $\hat k_1=r_1\geq k_1$ sections.
 By our choices of $\hat
k_i$,  $\sum_{i=1}^\ell \hat k_i\geq \sum_{i=1}^\ell k_i\geq k$.

 Let $0\leq\mu_1<\cdots<\mu_s\leq \mu <\mu_{s+1} <\cdots <\mu_\ell\leq 1$.
Suppose first that $\mu_s\neq 0$.
 By Lemma \ref{L:extension} there is a
nontrivial extension $0\to\widehat E_s\to G\to \widehat E_{s+1}\to 0$, with
$h^0(G)=\hat
k_s+\hat k_{s+1}$.  Thus,
$$
h^0\left(\widehat E_1\oplus\cdots\oplus\widehat E_{s-1}\oplus G\oplus \widehat
E_{s+1}\oplus \cdots\oplus \widehat E_\ell\right)=\sum_{i=1}^\ell \hat k_i\geq
k\ .
$$
On the other hand, by Proposition \ref{P:energyestimate} there is a hermitian
connection
on $\widehat E_1\oplus\cdots\oplus\widehat E_{s-1}\oplus G\oplus \widehat
E_{s+1}\oplus \cdots\oplus \widehat E_\ell$ and linearly independent sections
$\varphi_1,\ldots, \varphi_k$ such that $f_\tau(A,\vec\varphi)<
f_\tau(A_\infty, 0)=m$,
contradicting the minimality of $(A_\infty, 0)$.

Now suppose $\mu_s=\mu_1=0$, $\mu < \mu_i$ for $2\leq i\leq \ell$.  If for any
$2\leq i\leq
\ell$ we have $r_i < d_i+(r_i-\hat k_i)g$, then by Lemma \ref{L:extension}
there is a
nontrivial extension $0\to\widehat E_1\to G\to \widehat E_i\to 0$, with
$h^0(G)=\hat
k_1+\hat k_i$, and Proposition \ref{P:energyestimate} yields a contradiction as
before.
Suppose that for all $2\leq i\leq \ell$, $r_i=d_i+(r_i-\hat k_i)g$.  We claim
that
$\sum_{i=1}^\ell \hat k_i > k$.  For if $\sum_{i=1}^\ell \hat k_i = k$, then
$\sum_{i=2}^\ell (r_i-\hat k_i)= r-k$, and hence
$$
r > \sum_{i=2}^\ell r_i = \sum_{i=2}^\ell d_i+ (r_i-\hat k_i)g = d+ (r-k)g\ ;
$$
a contradiction.  Thus, we may replace $\widehat E_1$ by a bundle $\widehat
E_1^\prime$
having $\hat k_1^\prime=\hat k_1-1$ sections.  According to Lemma
\ref{L:extension}
 there is a
nontrivial extension $0\to\widehat E_1^\prime \to G\to \widehat E_2\to 0$, with
$h^0(G)=\hat
k_1^\prime +\hat k_2$, $\hat k_1^\prime +\sum_{i=2}^\ell\hat k_i\geq k$, and
Proposition
\ref{P:energyestimate} yields a contradiction as before.

\medskip
\noindent \emph{Case II}\@.  First notice that by the invariance of the
Yang-Mills-Higgs equations under the natural action by U($k$), we may assume
that
$\varphi_1^\infty, \ldots, \varphi_k^\infty$ form an $L^2$-orthogonal set of
sections.  In
particular, we may assume that there exists  $s\leq k$ such that
$\varphi_1^\infty,\ldots,\varphi_s^\infty$ are linearly independent and
$\varphi_{s+1}^\infty,\ldots,\varphi_k^\infty\equiv 0$.  Write
$E_\varphi=E_\varphi^\prime\oplus\Oka^{\oplus t}$, where $E_\varphi^\prime $
contains no
split factor of
$\Oka$.  Set $k_i=h^0(E_i)$,
$k_\varphi=h^0(E_\varphi)$, $k_\varphi^\prime
=h^0(E_\varphi^\prime)=k_\varphi-t$.  By
semicontinuity of cohomology, $k_\varphi+\sum_{i=1}^\ell k_i\geq k$.  As in
Case I, we may
replace each $E_i$ by a Hermitian-Yang-Mills bundle $\widehat E_i$ such that
$h^0(\widehat
E_i)=\hat k_i\geq k_i$, and \eqref{E:max} is satisfied for $i=1,\ldots,\ell$.
On the other
hand, since $E_\varphi$ satisfies the $k$-$\tau$-vortex equation for
$\tau=\mu+\gamma$ as
above, it follows that $E_\varphi^\prime$ is $\tau$-stable.  Therefore,
$0\neq \mu(E_\varphi^\prime)\leq \mu=\mu(E)$; and  since $\tau <
1$, we obtain from Lemma \ref{L:bounds} that $r_\varphi^\prime \leq
d_\varphi+(r_\varphi^\prime-k_\varphi^\prime)g$.  Finally, notice that since
$E^\infty$ is
unstable, $\mu_\ell > \mu$.  We may now apply Lemma \ref{L:extension} to
$E_\varphi^\prime$
and $\widehat E_\ell$ to obtain a nontrivial extension $0\to
E_\varphi^\prime\to G\to
\widehat E_\ell\to 0$, with $h^0(G)=k_\varphi^\prime+\hat k_\ell$.  It follows
that
$$
h^0\left(
G\oplus\Oka^{\oplus t}\oplus \widehat E_1\oplus\cdots\oplus\widehat E_{\ell-1}
\right) = k_\varphi +\sum_{i=1}^\ell \hat k_i \geq k\ .
$$
By Proposition \ref{P:energyestimate} there is a hermitian connection $A$ on
$G\oplus\Oka^{\oplus t}\oplus \widehat E_1\oplus\cdots\oplus\widehat
E_{\ell-1}$
and linearly independent sections $\varphi_1,\ldots,\varphi_k$ extending
$\varphi_1^\infty,\ldots,\varphi_s^\infty$ such that $f_\tau(A,\vec\varphi)<
f_\tau(A_\infty,\vec\varphi^\infty)=m$, again contradicting the minimality of
$m$.
This completes the proof of the Main Theorem.

\section{Stable Bundles}

We conclude by proving the Corollary stated in the Introduction.  Consider
first part
(ii).  The upper bound follows from Lemma \ref{L:bounds}.  By the Main Theorem,
it suffices
to show that if $E$ is semistable with $0<\mu < 1$ and $h^0(E)=r-1$, then $E$
is  stable.
Suppose to the contrary.  Then we can find a semistable subbundle $E'$ with
$Q=E/E'$ stable
and $\mu(E')=\mu(Q)=\mu$.  By Lemma \ref{L:bounds}, $h^0(E')\leq r'-1$, and
$h^0(Q)\leq r_Q-1$; contradiction.

Now consider part (i). Tensoring by ample line bundles allows us to restrict to
the case
$0<d\leq r$. Let
$\B$ be the set of gauge equivalence classes of solutions to the (one section)
$\tau$-vortex
equation for bundles of rank
$r$ and degree $d$.  By the proof of the Main Theorem, $\B\neq\emptyset$.  One
can therefore
show as in
\cite{BDW,BD} that
$\B$ is a smooth projective variety of dimension $(r^2-r)(g-1)+d$ with a
morphism $\psi :
\B\to\M$, where
$\M$ is the moduli space of semistable bundles of rank $r$ and degree $d$.  The
image of
$\psi$ is precisely the set of isomorphism classes of
 semistable bundles $E$ with $h^0(E)\geq 1$.
\begin{Lem}
Suppose that there exists a semistable (resp.\ stable) bundle $E_0$ of rank
$r$, degree
$d$, $0< d\leq r$, and  $h^0(E_0)=k\geq 1$.  Then there exists a semistable
(resp.\ stable)
bundle
$E$  of the same rank and degree with $h^0(E)=k-1$.
\end{Lem}

\begin{proof}  By Lemma \ref{L:bounds}, $k\leq r$.   The case where $d=r$ and
$E_0$ is
strictly semistable is trivial.  In the other cases, $k<r$, and we may  write
$$
\beta_0 : 0\to\Oka^{\oplus k}\to E_0\to F\to 0\ ,
$$
where by assumption the connecting homomorphism $\delta_0: H^0(F)\to
H^1(\Oka^{\oplus k})$
is injective.  Consider $\{L_t : t\in D\}$ a smooth local family of line
bundles
parametrized by the open unit disk $D\subset \CBbb$ and satisfying
$L_0=\Oka$ and $H^0(L_t)=0$, $t\neq 0$.  Set $G_t=\Oka^{k-1}\oplus L_t$.  The
semistability
of $E_0$ implies that $ H^0(F^\ast\otimes G_t)=0$. Hence, $\{ H^1(F^\ast\otimes
G_t) : t\in
D\}$ defines a smooth vector bundle $V$ over  $D$.  Let $\beta=\{\beta(t) :
t\in
D\}$ be a nowhere vanishing section of $V$ with $\beta(0)=\beta_0$.  Then
$\beta$ defines a
smooth family of nonsplit extensions $0\to G_t\to E_t\to F\to 0$ and a smooth
family of
connecting homomorphisms
$$
\delta_t : H^0(F)\lra H^1(G_t)\subset\Omega^{0,1}(U)\ ,
$$
where $U$ is the trivial rank $k$, $C^\infty$ vector bundle on $\Sigma$.  By
assumption,
$\delta_0$ is injective; hence, $\delta_t$ is injective for small $t$.  It
follows that
$h^0(E_t)=k-1$ for small $t$.  Furthermore, since $E_0$ is semistable (resp.\
stable) then
$E_t$ is also semistable (resp.\ stable) for small $t$.
\end{proof}

\noindent Since the condition $h^0(E)=1$  is open in $\B$, by the Lemma there
exists some
component
$\B^\prime$ of $\B$ containing an open dense set $\B^\ast$ consisting of pairs
$[E,\varphi]$
with
$h^0(E)=1$.  Let $W^\ast=\psi(\B^\ast)$.  We will assume that $W^\ast$ is
contained in the
strictly semistable locus of $\M$ and derive a contradiction.  We may assume
that each
irreducible component of $W^\ast$ is contained a subvariety $S$ parametrized by
bundles of
the form
$E_1\oplus E_2$, where $E_1$ is stable with $h^0(E_1)\geq 1$, $E_2$ is
semistable, and
$\mu(E_1)=\mu(E_2)=\mu(E)$.  It follows that
\begin{equation} \label{E:dim1}
\begin{aligned}
\dim S&=(r_1^2-r_1)(g-1)+d_1+r_2^2(g-1)+1 \\
 &=r^2(g-1)-2r_1r_2(g-1)-r_1(g-1)+d_1+1\ .
\end{aligned}
\end{equation}
For $[E]\in W^\ast\subset S$, we have
\begin{equation}
\label{E:dim2}
\dim_{[E,\varphi]}\B\leq \dim_{[E]} S +\dim \psi^{-1}([E]) \ .
\end{equation}
The dimension of the fiber of $\psi$ is given by
$h^1(E_2\otimes E_1^\ast)$.  Assume first that a generic $E_1$ is not
isomorphic to any
factor of
$E_2$. Then the fiber dimension  is $r_1 r_2(g-1)$.  Thus, we obtain from
\eqref{E:dim1} and \eqref{E:dim2} that
$$
(r^2-r)(g-1)+d \leq r^2(g-1)-r_1r_2(g-1)-r_1(g-1)+d_1+1\ ,
$$
or,
$$
r_2(r_1-1)(g-1)+ d_2-1 \leq 0\ .
$$
This yields a contradiction in all cases other than
$r=d=2$.  The latter situation is covered by the following construction which
is verified
by straightforward dimension counting:
\begin{Prop}
For generic line bundles $L$ of degree 2 and generic extensions $0\to\Oka\to
E\to L\to 0$,
$E$ is a stable rank 2 bundle of degree 2 with a nontrivial holomorphic
section.
\end{Prop}

\noindent  In case $E_1$ is isomorphic to some factor of $E_2$, the fiber
dimension
increases by 1.  On the other hand, in this case $W^\ast$ is contained in a
strict
subvariety of $S$, so by \eqref{E:dim2} the same argument applies.
Isomorphisms with more
factors are handled similarly. This completes the proof of the Corollary.

\end{document}